\documentclass[11pt
]{amsart}

\usepackage{hyperref}
\hypersetup{bookmarksdepth=3}

\usepackage{geometry}
\usepackage{amsmath,amssymb}
\usepackage{
mathrsfs}
\usepackage{esint}

\usepackage{tensor}

\usepackage{amssymb,amsmath,amsthm}

\newtheorem{theorem}{Theorem}

\theoremstyle{remark}

\newtheorem{proposition}[theorem]{Proposition}

\usepackage{graphicx}

\newcommand{\eq}[1]{\begin{align*}#1\end{align*}}  
\newcommand{\Span}{\operatorname{span}}

\def\BBR {{\mathbb R}}

\begin{document}

\title{Matrix weights, Littlewood Paley inequalities and the Riesz transforms}
\author{Nicholas Boros and Nikolaos Pattakos} 

\date{}

\begin{abstract}
{In the following, we will discuss weighted estimates for the squares of the Riesz transforms $R_{1}^2,..., R_{m}^{2}$ on $L^{2}(W)$ where $W\in\mathbb C^{d\times d}$ is an $A_{2}$ weight. We will show that if the ``Heat $A_{2}$ characteristic" of $W$ is sufficiently close to $1$ then there is a dimensional constant $c>0$ such that

$$\|R_{i}^{2}\|_{2,W}\leq 1+c\sqrt{[W]_{A^{h}_{2}}-1},$$
for all $i=1,...,m$. This is accomplished by proving a Littlewood--Paley estimate with the use of the Bellman function technique. These results follow the ones obtained in \cite{NV1}, \cite{NV2}, where similar questions were considered for the case of scalar Muckenhoupt $A_{p}$ weights $w$. The use of matrix weights turns out to be much more complicated than the scalar weights since the environment that we have to work in is no longer commutative.}
\end{abstract}

\maketitle


\begin{section}{introduction}
Throughout this article, we will denote $W$ as a positive, invertible and $L^{1}_{loc}(\mathbb R^m)$, $d\times d$ matrix with complex entries.  For such $W$ we define the heat extension to the upper half space $\mathbb R_{+}^{m+1}$ as the convolution of $W$ with the fundamental solution $k$ in $\mathbb R_{+}^{m+1}$ of the heat equation $u_{t}-\Delta u=0$,

$$W^{h}(x,t)=k_{t}\ast W(x)=\frac1{(4\pi t)^{\frac{m}{2}}}\int_{\mathbb R^m}W(y)\exp\Big(-\frac{|x-y|^{2}}{4t}\Big) \ dy.$$
We define the {\it heat $A_{2}$ characteristic} of $W$ to be the supremum

$$[W]_{A_{2}^{h}}=\sup_{(x,t)\in\mathbb R^{m+1}_{+}}\Big\|(W^{h}(x,t))^{\frac12}((W^{-1})^{h})^{\frac12}(x,t)\Big\|.$$
If this number is finite we say that $W\in A_{2}^{h}$. 

The components of the Riesz transform in $\BBR^m$, $R_{1},..., R_{m}$ are defined in the usual way, as the convolution operators

$$R_{i}f(x)=\frac{\Gamma(\frac{n+1}{2})}{\pi^{\frac{n+1}{2}}}p.v.\int_{\mathbb R^m}\frac{x_{i}-y_{i}}{|x-y|^{m+1}}f(y) \ dy,$$
$i=1, \dots ,m.$  It is known that for $i = 1, \dots m$ we have $\|R_{i}\|_{2}=1$, where by $\|...\|_{2}$ we denote the operator norm in the Hilbert space $L^{2}$. 

The main theorem of this paper is the following continuity result about matrix weights and the square of the Riesz transforms.

\begin{theorem}
\label{main}
There is a dimensional constant $c_{d}>0$ such that for all matrix weights $W$ of heat $A_{2}$ characteristic sufficiently close to $1$ the estimate
\begin{equation}
\label{Riesz}
\Big\|\sum_{i=1}^{n}\sigma_{j_{i}}R_{j_{i}}^{2}\Big\|_{2,W}\leq1+c_{d}\sqrt{[W]_{A_{2}^{h}}-1}
\end{equation}
is satisfied, where $\{ j_{1},j_{2}, \dots,j_{n}\}$ is an arbitrary subset of $\{1, \dots, m\}$ and $\sigma=\{\sigma_{j_{i}}\}_{i=1}^{n}$ is an arbitrary choice of signs.
\end{theorem}
Before we discuss the proof of Theorem \ref{main} let us make some comments. It has been proven by the last two authors that the same question of the continuity of the operator norm of a Calder\'on-Zygmund operator $T$ in $L^{p}(w)$ with respect to $w$, $1<p<\infty$, has a positive answer for scalar $A_{p}$ weights $w$ (for the Muckenhoupt $A_{p}$ classes see \cite{GCRF}). That is, if the scalar weight $w$ is ``close" to the constant weight $1$, then the number $\|T\|_{p,w}$ is close to $\|T\|_{p}$. See \cite{NV1} and \cite{NV2}. Notice that in our case the heat $A_{2}$ characteristic allows us to have a measure of determining when a matrix weight $W$ is ``close" to the constant weight $Id$ and that $[Id]_{A_{2}^{h}}=1$ since $\|k\|_{1}=1$. Hence, Theorem \ref{main} is a matrix analogue of the results appearing in \cite{NV1} and \cite{NV2}. Unfortunately, we only have the result for the square of the Riesz transforms and not for the Riesz transforms themselves. In the scalar case everything is much simpler since there is a very nice interplay between the $A_{p}$ classes and the $BMO$ space, which allows us to go almost back and forth between weights and functions of bounded mean oscillation. Such a relationship does not exist between matrix $A_{2}$ weights and $BMO$ because there is no ``good" analogue of $BMO$ for matrix valued functions. That is why we have to invent some ``new" path connecting closeness between $W$ and $Id$ and between the operator norm of an operator in $L^{2}(W)$ and $L^{2}$. The main ingredient of the proof is the following Littlewood-Paley type estimate which is obtained by using the Bellman function technique.

\begin{theorem}
\label{main1}
There is a dimensional constant $c_{d}>0$ such that for all matrix weights $W$ of heat $A_{2}$ characteristic sufficiently close to $1$ the estimate
\begin{equation}
\label{Little}
2\int_{\mathbb R^{m+1}_{+}}\sum_{i=1}^{m}\Big|\Big(\frac{\partial f^{h}(x,t)}{\partial x_{i}},\frac{\partial g^{h}(x,t)}{\partial x_{i}}\Big)_{\mathbb C^{d}}\Big| \ dx dt\leq(1+c_{d}\sqrt{[W]_{A_{2}^{h}}-1})\|f\|_{2,W}\|g\|_{2,W^{-1}}
\end{equation}
holds true for all vector functions $f,g\in C^{\infty}_{c}(\mathbb R^m)$ with values in $\mathbb C^d$.
\end{theorem}

It is known that in order to obtain an estimate as the one presented in Theorem \ref{main1} it suffices to find a suitable Bellman function. For this reason we also have the following theorem that guarantees that such a function exists.

\begin{theorem}
\label{main2}
For any $0<\delta<1$, define the domain $D_{\delta}=\{(X,Y,x,y,r,s)\in\mathbb R_{+}\times\mathbb R_{+}\times\mathbb C^{d}\times\mathbb C^{d}\times\mathbb C^{d\times d}\times\mathbb C^{d\times d}:
|(x,e)|\leq X^{\frac12}(se,e)^{\frac12}, |(y,e)|\leq Y^{\frac12}(re,e)^{\frac12}, 1\leq\|r^{\frac12}s^{\frac12}\|\leq 1+\delta\}$, for all directions $e\in\mathbb C^d$. Let $K$ be any compact subset of $D_{\delta}$. Then there exists a function $B=B_{\delta,K}(X,Y,x,y,r,s)$ which is infinitely differentiable in a small neighborhood of $K.$  Furthermore, for any $\epsilon>0$, $B_{\delta,K}$ can be chosen in such a way that the following estimates hold: 
\begin{align*}
& (A_{1}) \ \  0\leq B\leq(1+\epsilon)(1+c\sqrt{\delta})X^{1/2}Y^{1/2}\\
& (A_{2}) \ \  -d^{2}B\geq 2|(dx,dy)|,
\end{align*}
where $c$ is a constant that depends on the dimension $d$.
\end{theorem} 

This function will be obtained as a byproduct of a result first appeared in \cite{N}, which states that a similar statement to the one presented in Theorem \ref{main} holds for a version of the Martingale transform. More precisely, let us consider a matrix weight $W$ defined on $\mathbb R$ and for each dyadic subinterval $I$ of the real line the operator $<W>_{I}=\frac1{|I|}\int_{I}W$. It has a complete system of eigenvectors $e_{I}^{1},..., e_{I}^{d}$ which we can assume to be an orthonormal basis of $\mathbb C^d$ and let us also denote by 

$$h_{I}(x)=\frac{1}{\sqrt{|I|}}(\chi_{I_{+}}(x)-\chi_{I_{-}}(x)),$$
the Haar function associated to $I$ (we denote by $I_{+}$ the right ``son" of $I$ and by $I_{-}$ the left ``son"). Our Martingale $M_{\sigma}^{W}$, where $\sigma=\{\sigma_{I}\}$ is a collection of signs, acts on vectors functions $f$ in the following way

$$M_{\sigma}^{W}f(x)=\sum_{\substack{I\in\mathcal D \\1\leq k\leq d}}\sigma_{I}^{k}(f,h_{I}e_{I}^{k})_{2}\cdot h_{I}e_{I}^{k}.$$
In \cite{N} it was proven that if the $A_{2}$ characteristic of $W$, 

$$[W]_{A_{2}}=\sup_{I}\Big\|<W>_{I}^{\frac12}<W^{-1}>_{I}^{\frac12}\Big\|$$
is sufficiently close to $1$ then there is a dimensional constant $c_{d}>0$ such that we have the estimate
\begin{equation}
\label{eq}
\|M_{\sigma}^{W}\|_{2,W}\leq1+c_{d}\sqrt{[W]_{A_{2}}-1}.
\end{equation}
Now we have in our disposal all the ingredients of the proof which we are going to present in the next section. 

Let us elaborate on the motivation of studying weighted norm inequalities with matrix weights. The motivation of studying estimates of this type comes from stochastic processes and operator theory. Let us consider a multivariate random stationary process. For simplicity we consider the case of discrete time i.e. a sequence of $d$-tuples $x(n)=(x_{1}(n),...,x_{d}(n))$, $n\in\mathbb Z$, of scalar random variables such that $\mathbb E|x_{j}(n)|^{2}<+\infty$ and the correlation matrix 
$$Q(n,k)=\{Q(n,k)_{i,j}\}_{1\leq i,j\leq d}:=\{\mathbb Ex_{i}(n)\overline{x_{j}(n)}\}_{1\leq i,j\leq d},$$
depends only on the difference $n-k$ (we use the symbol $\mathbb E$ to denote the expectation). Without loss of generality we can assume that the process is complex valued. It is well known (see \cite{YR}) that there exists a matrix valued non-negative measure $M$ on the unit circle $\mathbb T$ whose Fourier coefficients coincide with the entries of the correlation matrix
$$Q(n,k)=\widehat{M}(n,k),$$
$n,k\in\mathbb Z$ and that if the process is completely regular then its spectral measure, $M$, is absolutely continuous with respect to the normalized Lebesgue measure $m$ on the unit circle, i.e. $dM=Wdm$. The past of the process is defined as
$$\mathcal X_{n}=\Span\{x_{j}(k):1\leq j\leq d, k<n\}$$
and the future as
$$\mathcal X^{n}=\Span\{x_{j}(k):1\leq j\leq d, k\geq n\}.$$
By writing ``span" we mean the closed linear span in the complex Hilbert space $L^{2}(\Omega,dP)$. If we consider the mapping
$$x_{j}(k)\mapsto z^{k}e_{j},$$
where $\{e_{j}\}_{1\leq j\leq d}$ is the standard orthonormal basis of $\mathbb C^d$, then we obtain an isometric isomorphism between $\Span\{x_{j}(k):1\leq j\leq d, k\in\mathbb Z\}$ and $L^{2}(W)$. The past and the future of the process are mapped to the subspaces of $L^{2}(W)$
$$X_{n}=\Span\{z^{k}\mathbb C^d:k<n\}$$
and
$$X^{n}=\Span\{z^{k}\mathbb C^d:k\geq n\},$$
respectively. In this representation the angle between past and future is nonzero if and only if the Riesz projection $P_{+}$ is bounded in the weighted space $L^{2}(W)$. All these applications and many more aspects of the theory of weighted inequalities with matrix weights are thoroughly discussed in the introduction of \cite{TV1} and the references therein. 

\end{section}


\begin{section}{the proofs of the theorems}

First we begin with the proof of Theorem \ref{main2} that claims the existence of a certain Bellman function.

\begin{proof}
Using duality in inequality (\ref{eq}) we get the following for all $f\in L^{2}(W)$ and $g\in L^{2}(W^{-1})$ and $J\in\mathcal D$

\eq{\frac1{4|J|}\sum_{I\in\mathcal D(J)}\sum_{k}|(<f>_{I_{+}},e_{I}^{k}) &-(<f>_{I_{-}},e_{I}^{k})||(<g>_{I_{+}},e_{I}^{k})-(<g>_{I_{-}},e_{I}^{k})||I|\\
& \leq (1+c\sqrt{[W]_{A_{2}}-1})\Big(\int_{J}(Wf,f)dx\Big)^{\frac12}\Big(\int_{J}(W^{-1}g,g) \ dx\Big)^{\frac12}.}
Let us take a look at the sum over $k$. For fixed $I$ the vectors $\{e_{I}^{k}\}_{k}$ constitute an orthonormal basis of $\mathbb C^{d}$ and the first summand of the sum over $k$ is the first coordinate of the vector $\Delta_{I}f:=<f>_{I_{+}}-<f>_{I_{-}}$ times the first coordinate of the vector $\Delta_{I}g$. The same thing happens until the last summand. Here we have to notice the following equality, which can then be estimated as

\eq{\sum_{k}|(<f>_{I_{+}},e_{I}^{k}) - (<f>_{I_{-}},e_{I}^{k})||(<g>_{I_{+}},e_{I}^{k}) &- (<g>_{I_{-}},e_{I}^{k})||I| \\
& =\sum_{k}|(\Delta_{I}f,e_{I}^{k})||(\Delta_{I}g,e_{I}^{k})||I|\\
& \ge |I|\Big|\sum_{k}(\Delta_{I}f,e_{I}^{k})(\Delta_{I}g,e_{I}^{k})\Big|\\ 
& = |I|\Big|(\Delta_{I}f,\Delta_{I}g)\Big|.}
So what we actually have proven is that 

\eq{\frac1{4|J|}\sum_{I\in\mathcal D(J)}\Big|\Big(<f>_{I_{+}} &-<f>_{I_{-}},<g>_{I_{+}}-<g>_{I_{-}}\Big)\Big||I|\\
& \leq (1+c\sqrt{[W]_{A_{2}}-1})\Big(\int_{J}(Wf,f) \ dx\Big)^{\frac12}\Big(\int_{J}(W^{-1}g,g) \ dx\Big)^{\frac12}.}
Hence, we can define the Bellman function to be the supremum of the left hand side

$$B(X,Y,x,y,r,s)=\sup\{\frac1{4|J|}\sum_{I\in\mathcal D(J)}\Big|\Big(<f>_{I_{+}}-<f>_{I_{-}},<g>_{I_{+}}-<g>_{I_{-}}\Big)\Big||I|:$$
$$x=<f>_{I}, y=<g>_{I}, X=\int_{J}(Wf,f) \ dx, Y=\int_{J}(W^{-1}g,g) \ dx, r=<W>_{J}, s=<W^{-1}>_{J}\}$$

over the domain

$$D_{Q}=\{(X,Y,x,y,r,s)\in\mathbb C^{d}\times\mathbb C^{d}\times\mathbb R_{+}\times\mathbb R_{+}\times\mathbb C^{d\times d}\times\mathbb C^{d\times d}:$$
$$|(x,e)|\leq X^{\frac12}(se,e)^{\frac12}, |(y,e)|\leq Y^{\frac12}(re,e)^{\frac12}, 1\leq\|r^{\frac12}s^{\frac12}\|\leq 1+\delta\},$$
for all directions $e$ in $\mathbb C^{d}$ and a fixed $\delta\in(0,1)$.

Obviously, the function $B$ satisfies 
\begin{equation}
\label{eq3}
0\leq B(X,Y,x,y,r,s)\leq (1+c\sqrt{\delta})X^{\frac12}Y^{\frac12}
\end{equation}

Notice that $B$ does not depend on the choice of the interval $J$. This happens because averages of functions remain the same under the use of a linear transformation. In addition, we claim that for all $6$-tuples $a^{+}=(X^{+}, Y^{+}, x^{+}, y^{+}, r^{+}, s^{+}), a^{-}=(X^{-},Y^{-},x^{-},y^{-},r^{-},s^{-})\in D_{\delta}$, such that $\frac{a^{+}+a^{-}}{2}\in D_{\delta}$, the following inequality is true
\begin{equation}
B\Big(\frac{a^{+}+a^{-}}{2}\Big)-\frac{B(a^{+})+B(a^{-})}{2}\geq\frac1{4}|(x^{+}-x^{-},y^{+}-y^{-})|.
\end{equation}
To prove this claim let us consider a positive $\epsilon$. There exists functions $f^{+}, g^{+}, W^{+}$ on $J_{+}$ such that they satisfy the conditions in the supremum of the function B for the vector $a^{+}$ and
$$B(a^{+})-\epsilon\leq\frac{1}{4|J_{+}|}\sum_{I\in\mathcal D(J_{+})}|(<f^{+}>_{I_+}-<f^{+}>_{I_-},<g^{+}>_{I_+}-<g^{+}>_{I_-})||I|.$$
Do the same for the vector $a^{-}$ in the interval $J_{-}$. Define the functions F, G, W on the interval J as: $F=\begin{cases} f_{+} &, \text{on}\, J_{+}\\
f_{-} &, \text{on}\, J_{-}\\
\end{cases},$
\qquad$G=\begin{cases} g_{+} &, \text{on}\, J_{+}\\
g_{-} &, \text{on}\, J_{-}\\
\end{cases},$
\qquad and,\qquad $W=\begin{cases} W_{+} &, \text{on}\, J_{+}\\
W_{-} &, \text{on}\, J_{-}\\
\end{cases}.$
Observe that they satisfy the required equalities in order to be acceptable for the supremum that defines the Bellman function for the vector $\frac{a^{+}+a^{-}}{2}$ and therefore,
\eq{B\left(\frac{a^{+}+a^{-}}{2}\right) &\geq\frac{1}{4|J|}\sum_{I\in\mathcal D(J)}|(F_{I+}-F_{I-},G_{I+}-G_{I-})||I|=\\
& \frac14\cdot\frac{1}{2|J_{+}|}\sum_{I\in\mathcal D(J_{+})}|(f^{+}_{I+}-f^{+}_{I-},g^{+}_{I+}-g^{+}_{I-})||I|\\
&+\frac14\cdot\frac{1}{2|J_{-}|}\sum_{I\in\mathcal D(J_{-})}|(f^{-}_{I+}-f^{-}_{I-},g^{-}_{I+}-g^{-}_{I-})||I|+
\frac{1}{4|J|}\text{term}(I=J)\\
& \geq\frac12(B(a^{+})-\epsilon)+\frac12(B(a^{-})-\epsilon)+\frac14|(x_{+}-x_{-},y_{+}-y_{-})|.}
Now we need to mollify this function $B$, in order to make it smooth. This can be done in exactly the same way as in \cite{NV}. The concavity inequality remains the same after the mollification and the size condition can become $1+C_{K}\epsilon$ times worse, where $C_{K}$ is a constant that depends on the compact set $K$. 
\end{proof}

The next step is to prove Theorem \ref{main1}. But, now that we have the Bellman function this proof is basically the same with the one given in \cite{PV}. The only difference is that we work with the function

$$v(x,t)=((Wf,f)^{h}(x,t),(W^{-1}g,g)^{h}(x,t),f^{h}(x,t),g^{h}(x,t),W^{h}(x,t),(W^{-1})^{h}(x,t)).$$
Finally, for the proof of Theorem \ref{main} we need to observe only the following result (for a proof see \cite{NV}).

\begin{proposition}
Let $\phi, \psi$ be real functions in $C_{c}^{\infty}(\mathbb R^m)$. Then the integral $\int \frac{\partial \phi^{h}}{\partial x_{i}}\cdot\frac{\partial \psi^{h}}{\partial x_{i}} \ dxdt$ converges absolutely for all $i=1, \dots,m$ and 
\eq{\int R_{i}^{2}\phi\cdot\psi dx=-2\int \frac{\partial \phi^{h}}{\partial x_{i}}\cdot\frac{\partial \psi^{h}}{\partial x_{i}} \ dx dt.}
\end{proposition}
Therefore, the left hand side of inequality (\ref{Riesz}) is bounded from above by the left hand side of inequality (\ref{Little}) and this means that we have the desired estimate.
\end{section}

\textbf{Acknowledgments}: The authors would like to thank professor Alexander Volberg from Michigan State University in East Lansing for useful discussions.

\vskip10pt

\noindent N. Boros:  Dept. of Mathematics, Olivet Nazarene University, Bourbonnais, Illinois, USA.\newline 
{\it e-mail address}:  \texttt{nboros@olivet.edu}

\vskip10pt

\noindent N. Pattakos:  School of Mathematics, University of Birmingham, Edgbaston, England.\newline
{\it e-mail address}:  \texttt{nikolaos.pattakos@gmail.com}

\end{document}